\documentclass[12pt,leqno]{article}

\usepackage{amssymb}
\usepackage{latexsym}

\newtheorem{theorem}{Theorem}[section]
\newtheorem{lemma}[theorem]{Lemma}

\newenvironment{proof}{\noindent {\bf
Proof}.\ }{\proofbox}

\def\numberlikeadb{\global\def\theequation{\thesection.\arabic{equation}}}
\numberlikeadb

\newcommand{\halmos}{\rule{1ex}{1.4ex}}
\newcommand{\proofbox}{\hspace*{\fill}\mbox{$\halmos$}}

\newcommand{\E}{\mathbb{E}}
\newcommand{\pr}{\mathbb{P}}
\def\re{\mathbb{R}}
\newcommand{\m}[1]{\marginpar{\tiny{#1}}}
\newcommand{\nats}{\mathbb{N}}

\newcommand{\dtv}{d_{\mbox{{\rm\tiny TV}}}}

\newcommand{\eqa}{\begin{eqnarray}}
\newcommand{\ena}{\end{eqnarray}}
\newcommand{\eq}{\begin{equation}}
\newcommand{\en}{\end{equation}}
\newcommand{\eqs}{\begin{eqnarray*}}
\newcommand{\ens}{\end{eqnarray*}}

\def\a{\alpha}
\def\b{\beta}
\def\g{\gamma}
\def\G{\Gamma}
\def\d{\delta}
\def\D{\Delta}
\def\e{\varepsilon}
\def\h{\eta}
\def\z{\zeta}
\def\th{\theta}

\def\m{\mu}
\def\n{\nu}
\def\p{\pi}
\def\r{\rho}
\def\s{\sigma}
\def\t{\tau}

\def\nin{\noindent}

\def\msk{\medskip}

\def\Blb{\left\{}
\def\Brb{\right\}}

\def\giv{\,|\,}

\def\non{\nonumber}

\def\Eq{\ =\ }
\def\Le{\ \le\ }

\def\Ref#1{{\rm (\ref{#1})}}

\def\bp{\begin{proof}}
\def\ep{\end{proof}}

\def\uk{^{(\KKK)}}

\def\ignore#1{}

\def\Xm{X^\m}

\def\tso{\t_{\{s,0\}}}

\def\tss{\t_{\{s\}}}

\def\law{{\cal L}}
\def\ep{\hfil\halmos\break\vskip20pt}

\def\mm{{\cal M}}

\def\jj{{\cal J}}

\def\sjj{\sum_{J\in\jj}}

\def\half{{\textstyle \frac12}}
\def\aaa{{\mathcal A}}

\def\remark{\nin {\bf Remark.\ }}

\def\qsd{quasi-stationary distribution}
\def\Def{\ :=\ }

\def\bir{b}
\def\dea{d}

\def\ui{^{(1)}}
\def\ut{^{(2)}}
\def\ex{\E}

\def\Bl{\left(}
\def\Br{\right)}
\def\alp{\alpha}
\def\BP{Barbour \& Pollett}
\def\pds{\p^{\d_s}}
\def\Xds{X^{\d_s}}

\def\KKK{\z}
\def\bbb{\b}
\def\ddd{\d}
\def\x{\xi}
\def\integ{\mathbb{Z}}
\def\S{\Sigma}
\def\trace{{\rm tr\,}}

\def\ts{{\tilde \s}}
\def\BBB{G}
\def\ttt{{\tilde\t}}
\def\htt{{\hat \t}}
\def\tpp{{\tilde\p}}

\def\uj{^{(j)}}
\def\tT{{\widetilde T}}
\def\tX{{\widetilde X}}
\def\tr{{\tilde r}}
\def\uds{^{\d_s}}
\def\tpds{\tpp\uds}

\begin{document}

\title{Total variation approximation for quasi-equilibrium distributions, II}
\author{
A. D. Barbour\footnote{Institut f\"ur Mathematik, Universit\"at Z\"urich,
Winterthurertrasse 190, CH-8057 Z\"URICH; and National University of Singapore.\msk}
\ and
P. K. Pollett\footnote{University of Queensland;
PKP was supported in part by the Australian Research Council Centre of Excellence
for Mathematics and Statistics of Complex Systems.
}\\
Universit\"at Z\"urich, National University of Singapore
\\ and University of Queensland
}

\date{}
\maketitle

\begin{abstract}
Quasi-stationary distributions, as discussed by Darroch \& Seneta (1965),
have been used in biology to describe the steady state behaviour of population
models which, while eventually certain to become extinct, nevertheless
maintain an apparent stochastic equilibrium for long periods.  These
distributions have some drawbacks: they need not exist, nor be unique,
and their calculation can present problems.
In an earlier paper, we gave biologically plausible conditions under which
the \qsd\ is unique, and can be closely approximated by distributions
that are simple to compute.  In this paper, we consider conditions under
which the \qsd, if it exists, need not be unique, but an apparent stochastic
equilibrium can nonetheless be identified and computed; we call such a
distribution a quasi-equilibrium distribution.
\end{abstract}

{\it Running head:}  Quasi-equilibrium distributions II

\section{Introduction}
\setcounter{equation}{0}
A rather general population growth model can be formulated as
a Markovian birth and death process~$X$ in continuous time,
where~$X(t)$ represents the number of individuals at time~$t$
in a population in a prescribed area~$A$, having transition rates
\eq\label{pop-model}
\begin{array}{cll}
    &q_{i,i+1} \Eq i\bbb(i/A),\qquad q_{i,i-1} \Eq i\ddd(i/A),\qquad\qquad &i\ge1;\\
    &q_{ij} \Eq 0 &\mbox{otherwise},
\end{array}
\en
where $\bbb(x)$ and~$\ddd(x)$ are the {\em per capita\/} rates of birth and
mortality at population density~$x=i/A$.
The logistic growth model of Verhulst~(1838) was the first to describe
mathematically the evolution of a population to a non-zero equilibrium,
contrasting with the Malthusian law of exponential growth, and
its stochastic version falls into the above framework, with $\bbb(x)=b$
constant in time, and with $\ddd(x) = d + cx$, linearly increasing in~$x$.
However, if the set of states $\nats$ is a communicating class and
$\inf_{x> 0} \ddd(x) > 0$, the stochastic model does not have a non-zero
equilibrium distribution even if $\bbb(0) > \ddd(0)$ and
$\lim_{x\to\infty} (\bbb(x)-\ddd(x)) < 0$, since $\nats$ is then transient,
and eventual absorption in~$0$ is certain.

Darroch and Seneta~(1965), building on the work of Yaglom~(1947) in
the context of branching processes, introduced the concept of a
quasi-stationary distribution, in an attempt to describe the long
term behaviour of a transient Markov chain prior to eventual
absorption. However, for chains with countably infinite state space,
Seneta and Vere--Jones~(1966) showed that the \qsd\ need neither
exist nor be unique.  Furthermore, even when there is a unique \qsd,
its calculation may pose substantial problems.  This apparently
makes the \qsd\ unsatisfactory for typical biological applications.
However, in \BP~(2010) [BP], we were able to give conditions, simply
expressed in terms of its properties, under which a continuous time
Markov chain~$X$ has exactly one \qsd. Under the same conditions,
the \qsd\ can be
approximated to a specified accuracy by the equilibrium distribution
of a positively recurrent `returned process'~$X^\m$, which may often
be much easier to compute. It was also shown, under slightly more
stringent conditions, that the distribution of~$X(t)$ is close to
its \qsd\ for long periods of time.

The conditions given in [BP] are satisfied for many population
models of the form~\Ref{pop-model}, including that of Verhulst~(1838).
However, a related model, in which the {\it per capita\/} death rate
$\ddd(x)=d$ remains constant as~$x$ increases, and the birth rate
declines exponentially, $\bbb(x) = be^{-\a x}$ for some $\a>0$ (Ricker 1954), does
not.  Indeed, although this biologically plausible model also gives rise to
apparently stable equilibrium behaviour for long periods of time, it
follows from van Doorn~(1991) that the process actually has infinitely many
possible \qsd s. To enable the long term behaviour of such models to be
adequately described, we now introduce
a new set of conditions, complementary to those in [BP], which can apply in
cases, such as that above, in which the \qsd\ need not exist nor be unique.

Denoting the state space of~$X$ by $C\cup\{0\}$, where~$0$ is the cemetery
state and~$C$ is irreducible, the returned process~$X^\m$ is also Markov.
It evolves like~$X$, except when it reaches the state~$0$.  Whenever it does,
instead of being absorbed in~$0$,
it is instantaneously returned to~$C$ according to the `return' probability
distribution~$\m$; hence each~$X^\m$ is a recurrent process.
Under our conditions, the returned processes for a wide class of return
distributions all have very similar equilibrium distributions,
and the distribution of~$X(t)$, given any reasonable fixed initial state, is also
similar to them for long periods of time.  Thus, for computational and practical
purposes, the situation is much as before.  The only difference is
that the \qsd\ can no longer be taken as the representative of the class of
`good' equilibrium distributions, since it need neither exist nor be unique.
Instead, any member~$\m$ of the class of `good' return distributions can be
chosen, and the equilibrium distribution of~$X^\m$ then serves as a
quasi-equilibrium
approximation to $\law(X(t))$ in the appropriate range of~$t$.

The main results, Theorems \ref{main-TV} and~\ref{quasi-longtime}, are proved in
Section~\ref{return-measure}.    In Section~\ref{b&d}, as an illustration, we discuss
the application of the theorems to birth and death processes.  These processes have
been widely studied, because of their relatively simple structure, and allow
our results to be easily interpreted.  Our theorems are however equally applicable
to processes with more complicated structure, and we illustrate their application
to Markov population processes in several dimensions in Section~\ref{MPP}.

\section{The return approximation}\label{return-measure}
\setcounter{equation}{0}
Assume that~$X$ is a stable, conservative and non-explosive pure jump Markov process
on a countable state space, consisting of a single transient class~$C$ together
with a cemetery state~$0$. For any probability distribution~$\m$ on~$C$, define
the modified process~$\Xm$ with state space~$C$ to have exactly the
same behaviour as~$X$ while in~$C$, but, on reaching~$0$, to be instantly returned
to~$C$ according to the distribution~$\m$.  Thus, if~$Q$ denotes the infinitesimal matrix
associated with~$X$, and~$Q^\m$ that belonging to~$X^\m$, we have
\eq\label{Qmu-def}
   q_{ij}^\m \Eq q_{ij} + q_{i0}\m_j \qquad \mbox{for } i,j \in C.
\en
In this section, under a rather simple set of conditions, we show that
the {\em stationary\/} distribution $\p^\m$ of~$X^\m$ is little influenced by the
choice of~$\m$, for~$\m$ in a large class~$\mm$ of distributions.
We give a bound, uniform for all $\m,\n \in\mm$, on the total variation distance
$$
    \dtv(\p^\n,\p^\m) \Def \sup_{A\in C}|\p^\n\{A\} - \p^\m\{A\}|
        \Eq \frac12 \sum_{k\in C}|\p^\n(k) - \p^\m(k)|
$$
between $\p^\n$ and~$\p^\m$, that is expressed in terms
of hitting probabilities and mean hitting times for the process~$X$.
The bound is such that it can be expected to be small in circumstances
in which the process~$X$ typically spends a long time in~$C$ in apparent equilibrium,
before being absorbed in~$0$ as a result of an `exceptional' event.

Define
\eq\label{tau-def}
  \t_A \ :=\ \inf\{t>0\colon\,X(t) \in A, \,X(s) \notin A \mbox{ for some } s < t\},
\en
with the infimum over the empty set being taken to be~$\infty$, noting
that $\t_A > 0$ a.s.\ even when $X(0) \in A$.
Our basic conditions can then be expressed as follows.

\medskip
\nin{\bf Condition B}.\ \,{\sl There exists $s\in C$ such that, defining
\eqs
    p_k &:=& \pr_k[X_{\tso} = s];\quad
    T_k \Def \ex_k[\tso],
\ens
we have}
\eqs
     &{\rm (i)}&  \inf_{k\in C} p_k  \Eq p > 0\,;\\
     &{\rm (ii)}&  T_k  \ <\ \infty \ \mbox{ for all } k \in C .
\ens

\nin Here, $\pr_k$ and $\E_k$ refer to the distribution of~$X$
conditional on $X(0)=k$.

  Condition B(i)
can be expected to be satisfied in reasonable generality, and is the same as
Condition A(i) in [BP]. Condition B(ii) substantially relaxes Condition A(ii)
in [BP], which stipulated that $T_k \le T < \infty$, {\it uniformly\/}
for all~$k\in C$.  If~$X$ is typically to spend a long time in apparent
equilibrium before being absorbed in~$0$, it will be necessary for $1-p_s$,
the probability that an excursion from~$s$ lands in~$0$, to be small.

We first note that
\[
  T_s \Eq \sum_{k\in C} T_{sk} \ <\ \infty,
\]
where
\[
  T_{sk} \Def \int_0^\infty \pr_s[\{\tso > t\}\cap \{X(t)=k\}] \,dt
\]
is the expected amount of time spent in~$k$ before first returning to $\{s,0\}$,
starting in~$s$.  Hence, for any $\KKK > 0$, we can pick $C_\KKK \subset C$
such that
\eq\label{CK-def}
   \sum_{k\notin C_\KKK} T_{sk} \Le \KKK(1-p_s)T_s;
\en
we do so in such a way that $s\in C_\KKK$, and that $T^+_\KKK := \sup_{k\in C_\KKK} T_k$
is as small as possible. We then define the process~$X_\KKK$ to be
the same as~$X$, except that any excursions outside~$C_\KKK$ take zero time to
complete. This process~$X_\KKK$ now satisfies Condition~A of~[BP], so that the
results of~[BP] can be applied to it.  Finally, we extend
the results for~$X_\KKK$ to the process~$X$.  To accomplish this programme, we need some
preparatory results.

\begin{lemma}\label{new-2.2}
Define $\m(T) := \sum_{k\in C} \m(k)T_k$. Then, under Condition~B,
\eqs
   &{\rm (i)}& \m(T) \Le \ex^\m \tss^\m \Le \m(T)/p\,;\\
   &{\rm (ii)}&  \ex_k \tss^\m \Le T_k + (1-p_k)\m(T)/p\,, \quad k \in C;\\
   &{\rm (iii)}& \ex_s \tss^\m \ <\ \infty\ \mbox{ if and only if } \m(T) < \infty\,,
\ens
where $\t_A^\m$ is defined similarly to~$\t_A$, but with the process~$\Xm$ in place of~$X$,
and $\ex^\m$ denotes expectation under the initial distribution~$\m$.
\end{lemma}

\proof
The proof is based on the equation
\eq\label{hitting-times}
   \tss^\m \Eq \tso^\m + \tss^{\m,1}, 
\en
in which $\tss^{\m,1}$ is the time that elapses after~$\tso^\m$ until~$\Xm$ first
reaches~$s$, zero if $\Xm(\tso) = s$.  Taking expectations with respect to~$\pr^\m$,
this yields
\[
   \ex^\m \tss^\m \Eq \sum_{k\in C} \m(k)\{T_k + (1-p_k)\ex^\m \tss^{\m,1}\},  
\]
from which it follows that
\[
    \m(T) \Le \ex^\m \tss^\m \Le \m(T) + (1-p)\ex^\m \tss^\m
\]
and Part~(i) is proved.  Part~(ii) follows by taking expectations in~\Ref{hitting-times}
with respect to~$\pr_k$, which also gives
\[
   \ex_k\tss^\m \ \ge\ (1-p_k)\ex^\m \tss^\m \ \ge\ (1-p_k)\m(T).
\]
Part~(iii) follows from these considerations, taking $k=s$.
\ep

We now define
\[
   \mm_M \Def \{\m \in {\cal PM}(C)\colon\, \m(T) \le MT_s\},
\]
for any $M>0$.  The next lemma bounds the equilibrium probability
that $\Xm \notin C_\KKK$, for any $\m\in\mm_M$.

\begin{lemma}\label{outside-CK}
  Under Condition~B, for any $\m \in \mm_M$, we have
\[
   \p^\m(C_\KKK^c) \Le (1-p_s)\{\KKK + M/p\} \ =:\ \e(\KKK,M).
\]
\end{lemma}

\proof
By a standard renewal argument,
\eqs
   \p^\m(A)\ex_s\tss^\m &=& \sum_{k\in A} \ex_s
      \Blb \int_0^\infty I[\tss^\m > t]\, I[X(t)=k]\,dt \Brb \\
   &\le& \sum_{k\in A} T_{sk} + (1-p_s)\ex^\m \tss^\m .
\ens
It thus follows from~\Ref{CK-def}, \Ref{hitting-times} and Lemma~\ref{new-2.2}(i)
that
\[
  \p^\m(C_\KKK^c)T_s \Le \p^\m(C_\KKK^c) \ex_s\tss^\m \Le (1-p_s)(\KKK T_s + \m(T)/p),
\]
and the lemma follows.
\ep

In what follows, we assume that
$M\ge1$, ensuring that the distribution~$\d_s$ that puts
probability~$1$ on the state~$s$ itself belongs to~$\mm_M$.

We now return to the pure jump Markov process $X_\KKK$, which has the same jump
chain as~$X$, and the same jump rates~$q_k$ for all $k\in C_\KKK$, but with
$q_k=\infty$ for $k\notin C_\KKK$.  We also define its returned processes~$X_\KKK^\m$
in the same way as for~$X$, but with the new jump rates.  We then define
\eq\label{Tsuk-def}
   T_s\uk \Def \sum_{k\in C_\KKK} T_{sk} \ \ge\ T_s\{1 - \KKK(1-p_s)\},
\en
the mean time for~$X_\KKK$ to return to the set~$\{0,s\}$,
starting from~$s$, the last inequality following from~\Ref{CK-def}.

\begin{theorem}\label{main-TV}
   Suppose that Condition~B holds, and that $M\ge1$.
 Then, for any $\m\in\mm_M$,
\[
   \dtv(\p^\m,\pds) \Le 2(1-p_s)\Bl \frac{T^+_\KKK}{pT_s} + \KKK + \frac M{p} \Br.
\]
\end{theorem}

\proof
We begin by considering the process $X_\KKK^\m$ for any $\m\in\mm_M$, noting
that, for any $k\in C_\KKK$, its equilibrium distribution~$\p_\KKK^\m$ satisfies
\eq\label{K-measure}
    \p_\KKK^\m(k) \Eq \p^\m(k)/\p^\m(C_\KKK).
\en
Now the
process~$X_\KKK$ satisfies Condition~A of~[BP], and hence, from (2.13) of~[BP],
\eq\label{BP-thm}
   \dtv(\p_\KKK^\m , \pds_\KKK) \Le 2(T^+_\KKK/p)\sum_{k\in C_\KKK}\pds_\KKK(k)q_{k0}.
\en
Then, by a renewal argument, letting $N_k(t)$ denote the number of
visits of~$X^{\d_s}$ to~$k$ in $[0,t]$, we have
\eqa
  \sum_{k\in C}\pds(k)q_{k0} &=& \lim_{t\to\infty} t^{-1}N_0(t) \non\\
    &=& \lim_{t\to\infty} \{t^{-1}N_s(t)\} \lim_{t\to\infty} \{N_0(t)/N_s(t)\}
    \Eq T_s^{-1} (1-p_s). \label{ergodic}
\ena
It now follows from \Ref{K-measure}--\Ref{ergodic} that
\[
    \dtv(\p_\KKK^\m , \pds_\KKK) \Le 2(T^+_\KKK/p) T_s^{-1}(1-p_s)/\pds(C_\KKK).
\]
Hence
\eqs
   \lefteqn{\dtv(\p^\m , \pds)}\\
   &=& \frac12 \sum_{k\in C_\KKK}|\p_\KKK^\m(k)\p^\m(C_\KKK) - \pds_\KKK(k)\pds(C_\KKK)|
     + \frac12 \sum_{k\notin C_\KKK}|\p^\m(k) - \pds(k)| \\
   &\le& \pds(C_\KKK)\dtv(\p_\KKK^\m , \pds_\KKK) + \half|\pds(C_\KKK) - \p^\m(C_\KKK)|
       + \half(\p^\m(C_\KKK^c) + \pds(C_\KKK^c)) \\
   &\le& 2(T^+_\KKK/p) T_s^{-1}(1-p_s) + 2\e(\KKK,M),
\ens
this last from Lemma~\ref{outside-CK}; as before, $\d_s\in\mm_M$,
because $M\ge1$.
\ep

\remark Of course, for the theorem to imply that~$\p^\m$ and~$\pds$ are close,
one needs~$(1-p_s)$ to be very small, which has already been noted as a necessary
condition for long time stability.  One also needs
$\frac{T^+_\KKK}{pT_s} + \KKK + \frac M{p}$
not to be too large.  The smaller~$\KKK$ is chosen, the larger is the value of~$T^+_\KKK$,
so that, in specific models, there is an optimum choice of~$\KKK$, limiting the
accuracy of approximation that can be demonstrated by this method.

\medskip
We now turn our attention to the distribution of~$X(t)$ for fixed values of~$t$,
starting from any particular state in~$C_\KKK$, and compare it to~$\pds$.
We begin by taking the initial state of~$X$ to be~$s$, and remark later that this
restriction makes little difference, provided that~$s$ is hit at least once.
To state the theorem, we define
\[
   r_\KKK \Def \pr_s[X \mbox{ does not leave } C_\KKK \mbox{ or hit } 0
        \mbox{ before returning to } s].
\]
Since $\lim_{\KKK\to0} r_\KKK = p_s$, the quantity $1-r_\KKK$ can be made as close
as desired to $1-p_s$ by decreasing~$\KKK$, but at the cost of increasing~$T_\KKK^+$
at the same time.  A crude bound for $1-r_\KKK$ in terms of $1-p_s$ comes from
observing that
\[
   \sum_{k\notin C_\KKK} T_{sk} \ \ge\ (1-r_\KKK)/q\uk,
\]
where $q\uk := \sup_{k\in J_\KKK}q_k$ and
\[
J_\KKK := \{k\notin C_\KKK\colon\, q_{kj} > 0 \mbox{ for some } j \in C_\KKK\};
\]
from~\Ref{CK-def}, this gives
\[
    (1-r_\KKK) \Le \KKK q\uk T_s(1-p_s).
\]

\begin{theorem}\label{quasi-longtime}
  Suppose that Condition~B holds, and let $B_\KKK := T_\KKK^+ q_s/p$.
\ignore{
\smallskip
\nin (1) There is a universal constant~$D$ such that
\[
  \dtv(\pds_\KKK,\law_s(X(t))) \Le (1-r_\KKK)(t/T_s\uk)
       + DB_\KKK\sqrt{ \frac {T_\KKK^+}{pt}} + (2/e)^{pt/16T_\KKK^+},
\]
for all $t\ge 16T_\KKK^+/p$.
\smallskip
\nin (2)
}
If $\e(\KKK,1) \le 1/2$, then, for all $t\ge 16T_\KKK^+/p$,
\eqs
  \lefteqn{\dtv(\law_s(X(t)),\pds)}\\ &\le&
     (1-r_\KKK)(2t/T_s + \KKK + 1/p) + DB_\KKK\sqrt{ \frac {T_\KKK^+}{pt}} + (2/e)^{pt/16T_\KKK^+}
    \ =:\ \h_\KKK(t).
\ens
\end{theorem}

\remark
Hence, informally, if $(1-r_\KKK)B_\KKK^2T_\KKK^+/(pT_s) \ll 1$ and $(1-r_\KKK) \ll 1$,
the distribution $\law_s(X(t))$ is close to~$\pds$ for all times~$t$ such that
\[
    B_\KKK^2 T_\KKK^+/p\ \ll\ t\ \ll\ T_s/(1-r_\KKK);
\]
note that $B_\KKK\ge1$, so that then $t \gg T_\KKK^+/p$ also.

\medskip

\proof
The argument is based on coupling two copies $X_\KKK\ui$ and~$X_\KKK\ut$ of the return
process~$\Xds_\KKK$, with~$X_\KKK\ui$ in equilibrium and with~$X_\KKK\ut$ starting in~$s$,
by the method used in [BP], Theorem~2.5.
The coupling is achieved by forcing~$X_\KKK\ui$ to follow the same sequence of states
as~$X_\KKK\ut$ after the first time that it hits~$s$, and to have identical residence
times in all states other than~$s$; the careful matching of the exponentially distributed
residence times of the two processes in~$s$ is all that is used to achieve the
coupling. Now the argument leading to (2.18) of [BP] shows that $X_\KKK\ui$
and~$X_\KKK\ut$ can be jointly
defined in such a way that, if $t\ge 16T_\KKK^+/p$, the event
$\D_{\KKK t}$ that they have coupled by~$t$ is such that
\[
  \pr[\D_{\KKK t}^c] \Le
          4c_G B_\KKK\sqrt{ \frac {T_\KKK^+}{pt}} + (2/e)^{pt/16T_\KKK^+},
\]
for a universal constant~$c_G$, not depending on~$\KKK$.  Now, because~$X_\KKK\ui$
is in equilibrium,
\[
  \pr[X_\KKK\ui\ \mbox{hits}\ \{0\}\cup C_\KKK^c\ \mbox{before}\ t]
    \Le t \sum_{k\in C_\KKK} \pds_\KKK(k)\sum_{l\in \{0\}\cup C_\KKK^c} q_{kl},
\]
and the double sum is bounded by $(1-r_\KKK)/T_s\uk$, as in the argument
leading to~\Ref{ergodic}.  If~$X_\KKK\ui$ does not hit $\{0\}\cup C_\KKK^c$ before~$t$,
and if $\D_{\KKK t}$ holds, then~$X_\KKK\ut$ also avoids $\{0\}\cup C_\KKK^c$
up to time~$t$, in which case it is indistinguishable from an $X$-process,
starting in~$s$. It thus follows that
\eq\label{truncated-bnd}
   \phantom{HH}\dtv(\pds_\KKK,\law_s(X(t))) \Le (1-r_\KKK)(t/T_s\uk)
       + DB_\KKK\sqrt{ \frac {T_\KKK^+}{pt}} + (2/e)^{pt/16T_\KKK^+},
\en
with $D := 4c_G$.  To complete the proof,
it now merely remains to note that $\dtv(\pds_\KKK,\pds) = \pds(C_\KKK^c)$,
and to use Lemma~\ref{outside-CK}; note that $(1-p_s) \ge (1 - r_\KKK)$,
and that $T_s\uk \ge T_s/2$, from~\Ref{Tsuk-def}, if $\e(\KKK,1) \le 1/2$.
\ep

\remark
Denoting by $A_\KKK$ the event that~$X$ hits $s$ before~$\{0\}\cup C_\KKK^c$,
the same argument can be used to show that
$\dtv(\law_k(X(t) \giv A_\KKK),\law_s(X(t)))$ is at most~$\h_\KKK(t)$
for any $k\in C_\KKK$, under the conditions of Theorem~\ref{quasi-longtime}.
Hence, conditional on the event that $X$ hits~$s$
before reaching~$\{0\}\cup C_\KKK^c$,  the distribution of~$X(t)$ starting from any $k\in C_\KKK$
is also close to~$\pds$ for all times~$t$ such that
\[
    B_\KKK^2 T_\KKK^+/p\ \ll\ t\ \ll\ T_s/(1-r_\KKK),
\]
provided that $(1-r_\KKK)B_\KKK^2T_\KKK^+/(pT_s) \ll 1$.  Thus the return
distribution~$\pds$ is then indeed an
appropriate long time approximation to the distribution of~$X$ in~$C$,
for times $t \ll T_s/(1-r_\KKK)$, and~$\pds$ can
be replaced by $\p^\m$ for any~$\m$ such that $\m(T) < \infty$, with extra
error at most that given by the bound in Theorem~\ref{main-TV},
with $\m(T)/T_s$ for~$M$.

The emphasis up to now has been on approximating $\law(X(t))$ by~$\pds$.
However, there are times when this approximation may also not be useful.
Examples of this are processes in which a set~$C_\KKK$ can be found that
has the properties that $T_s\uk$ and $T_\KKK^+$ are only moderately large
and $(1-r_\KKK)$ is tiny, but for which~\Ref{CK-def} is {\em not} satisfied.
Such is the case if there are states~$k \notin C_\KKK$ such
that~$T_k$ is extremely large; for instance,
if the equilibrium around~$s$ is metastable, $T_s$ itself may be enormously
larger than~$T_s\uk$.  Here, nonetheless, the intermediate bound~\Ref{truncated-bnd}
shows that~$\pds_\KKK$ acts as a good approximation for very long periods,
even though~$\pds$ may be very different.

In practice, computing~$\pds_\KKK$
may be complicated by having to cope with the detail of the return distribution
from~$C_\KKK^c$, which should not really be relevant here.  The final
approximation
is therefore phrased instead in terms of the accelerated return
process~$\tX\uds_{C'}$, for some $C'\subset C$ containing~$s$ but not~$0$, which
is returned to~$s$
at each time of leaving~$C'$.  Here, the set~$C'$ may reasonably be chosen to be finite,
in which case computing the equilibrium distribution~$\tpds_{C'}$ of the accelerated return process
becomes relatively easy.  We now define
 $\tT_{k,C'} := \ex_k[\ttt\uds_{\{s\}}]$, where~$\ttt\uds$ is defined as
in~\Ref{tau-def}, but with the process~$\tX\uds_{C'}$ in place of~$X$; and we set
$\tT_{C'}^+ := \sup_{k\in C'} \tT_{k,C'}$ and
$\tr_{C'} := \pr_s[\ttt\uds_{C\setminus C'} = \ttt\uds_{\{s\}}]$.

\begin{theorem}\label{quasi-longtime-2}
  Suppose that Condition~B\,(ii) holds, and let $B_{C'} := \tT_{C'}^+ q_s$.
Then
\[
  \dtv(\tpds_{C'},\law_s(X(t))) \Le (1-\tr_{C'})(t/\tT_{s,C'})
       + DB_{C'}\sqrt{ \frac {\tT_{C'}^+}{t}} + (2/e)^{t/16\tT_{C'}^+},
\]
for all $t\ge 16\tT_{C'}^+$, with~$D$ the same constant as in Theorem~\ref{quasi-longtime}.
\end{theorem}

\proof
The argument runs exactly as in the proof of~\Ref{truncated-bnd}, but with the
process~$\tX\uds_{C'}$ instead of $X_\KKK\uds$.  Since this process has no
absorbing state~$0$, $p$ can be replaced by~$1$ in the bound.
\ep

Theorem~\ref{quasi-longtime-2} is very much in line with the main message
of the paper.  The difference between Condition~A of~[BP] and Condition~B of
this paper largely concerns properties of the process starting from states
that it rarely ever reaches, and such differences should not prevent effective
approximation of the distribution of the process, at least for long periods
of time.  The essential difference between the situation in which Condition~A
is satisfied and that in which it is not is that, when it is not satisfied,
the approximating distribution need not be a quasi-stationary distribution
of the process, or even one of its return distributions, but instead a return
distribution associated with the process restricted to a truncated state space.
We consider an example of this in Section~\ref{MPP}.

\section{Birth and death processes}\label{b&d}
\setcounter{equation}{0}

Let~$X$ be a birth and death process with
birth rates $\bir_j \ge 0$, $1\le j < \infty$, with $\bir_0 = 0$, and with strictly
positive death rates $\dea_j$,
$j\ge0$. Define $\alp_1=1$ and
$$
\alp_j \Eq \frac{\bir_1 \cdots \bir_{j-1}}{\dea_2 \cdots \dea_j},\qquad j\ge1;
$$
then set
$$
  S_r^m \Def \sum_{l=r}^m \frac1{\alp_l\dea_l}.
$$

In order to use the theorems of the previous section, we need to find expressions
for the quantities $p$, $p_s$, $T_{sk}$, $T_s$, $T_\KKK^+$ and~$r_\KKK$ that appear
there.  These can be derived using hitting probabilities, which can be simply
expressed using the $\alp_j$ and the~$S_r^m$.  First, for any $j < m < l$ ,
we have
\eq\label{basic-hit}
   \pr_m[X\ \mbox{hits}\ l\ \mbox{before}\ j] \Eq S_{j+1}^m / S_{j+1}^l.
\en
A first consequence is that
\eqa
   1 - p_s &=& \frac{\dea_s}{\bir_s+\dea_s}\Bl 1 - \frac{S_1^{s-1}}{S_1^s} \Br
      \Eq \frac{1}{\alp_s(\bir_s+\dea_s)S_1^s}; \label{ps} \\
    p &=& p_1 \Eq \frac1{d_1S_1^s}.\label{pmin}
\ena
Next, if $i \notin \{0,s\}$, write $u_{ki} \Def \pr_k[\t_{\{i\}} < \tso]$,
$k\ne i$, and $u_{ii} = 1$: then we have
\eqa
    u_{ki} &=& \left\{
     \begin{array}{ll}
                0   &\quad\mbox{if}\ k < s < i\ \mbox{or}\ i < s < k;\\
                1   &\quad\mbox{if}\ s < i \le k;\\
                S_{s+1}^k/S_{s+1}^i   &\quad\mbox{if}\ s < k \le i;\\
                S_{1}^k/S_{1}^i   &\quad\mbox{if}\ 0 < k \le i < s;\\
                S_{k+1}^s/S_{i+1}^s   &\quad\mbox{if}\ 0 < i \le k < s,
               \end{array}\right. \label{uki}
\ena
from which it follows that, for such~$i$,
\eq\label{return-prob}
   1 - \pr_i[\t_{\{i\}} < \tso] \Eq \left\{
     \begin{array}{cl}
            \frac1{\alp_i(\bir_i+\dea_i)}\Blb \frac1{S_{i+1}^s} + \frac1{S_1^i} \Brb
               &\quad\mbox{if}\ 0 < i < s;\\ [2ex]
            \frac1{\alp_i(\bir_i+\dea_i)S_{s+1}^i}  &\quad\mbox{if}\ i > s.
        \end{array}\right.
\en
Also, for $i \notin \{s,0\}$ and $k \ne s$, we have
\eqa
   T_{ki} &:=& \int_0^\infty \pr_k[\{\tso > t\}\cap \{X(t)=i\}] \,dt \non\\
    &=& \frac{u_{ki}}{\bir_i+\dea_i}\,\frac1{1 - \pr_i[\t_{\{i\}} < \tso]},\label{Tki}
\ena
with $T_{ks} = T_{k0} = 0$, and then $T_{ss} = 1/(\bir_s+\dea_s)$ and
\eq\label{Tsi}
   T_{si} \Eq \left\{
     \begin{array}{cl}
        \frac{\dea_s T_{s-1,i}}{\bir_s+\dea_s}  &\quad\mbox{if}\ 0 < i < s;\\[2ex]
        \frac{\bir_s T_{s+1,i}}{\bir_s+\dea_s}  &\quad\mbox{if}\  i > s.
     \end{array}\right.
\en
These in turn give
\eq\label{Tk}
  T_k \Eq \sum_{i\ge1}T_{ki} \Eq \left\{
     \begin{array}{cl}
        \frac1{S_1^s}\sum_{i=1}^{s-1} \alp_i S_1^{i\wedge k} S_{(i\vee k)+1}^s
            &\quad\mbox{if}\ 0 < k < s;\\[2ex]
       \sum_{i\ge s+1} \alp_i S_{s+1}^{i\wedge k} &\quad\mbox{if}\ k > s,
    \end{array}\right.
\en
and
\eq\label{Ts} \mbox{}\qquad
   T_s \Eq (\bir_s T_{s+1} + \dea_s T_{s-1} + 1)/(\bir_s+\dea_s) \Eq
     \frac1{\alp_s(\bir_s+\dea_s)} \sum_{i \ge1} \alp_i S_1^{i\wedge s}/S_1^s.
\en
Now, choosing any value of $s > 0$, the formulae \Ref{Tsi}, \Ref{ps} and~\Ref{Ts}
can be used for any~$\KKK$ to determine a suitable set $C_\KKK := \{1,2,\ldots,a_\KKK\}$
so that~\Ref{CK-def} is satisfied, and~\Ref{Tk} enables both $T_\KKK^+$ and~$\m(T)$
to be determined. Furthermore, it follows from~\Ref{return-prob} with
$a_\KKK$ for~$s$ and with $s$ for~$i$ that
\eq\label{rK}
   1 - r_\KKK \Eq
     \frac1{\alp_s(\bir_s+\dea_s)}\Blb \frac1{S_{s+1}^{a_\KKK}} + \frac1{S_1^s} \Brb.
\en
Thus, and from~\Ref{pmin}, all the ingredients for the bounds
in Theorems \ref{main-TV} and~\ref{quasi-longtime} are available, recalling
also, for the calculation of~$B_\KKK$, that $q_s = \bir_s + \dea_s$.

For example, take the birth and death process given in~\Ref{pop-model}
with~$A$ large, $\ddd(x) = d$ constant, and with~$\bbb(\cdot)$ given by the Ricker
choice $\bbb(x) = be^{-\a x}$; thus $\bir_j = jbe^{-\a j/A}$ and $\dea_j = jd$.
If $b>d$, the deterministic equilibrium, in which the birth and death rates are
equal, is given by
$x = \frac1{\a}\log(b/d) =: c > 0$, suggesting the choice of $s :=
s(A) := \lfloor Ac \rfloor$.  This gives
\eqa
   \alp_j \ \sim\ j^{-1}(b/d)^j e^{-\a j(j+1)/2A}; \qquad \alp_{s} \ \asymp\ s^{-1}(b/d)^{s/2};\non\\
   \alp_j/\{s\alp_s\} \ \asymp \ j^{-1}e^{-\a j/2A} e^{-\a (j-s)^2/2A}; \label{calcs}
\ena
note that $\alp_s$ is exponentially large in~$A$. Thus,
immediately,  $S_1^s \ge 1/d$, and $S_{s+1}^a \ge 1/d$ if
 $a >  2(s(A)+1)$.  Hence $1 - p_s  = O\{(b/d)^{-s/2}\}$ from~\Ref{ps},
and $1-r_\KKK =O\{ (b/d)^{-s/2}\}$ also if $a_\KKK \ge 2(s(A)+1)$, from~\Ref{rK}.
Furthermore, $S_1^s$ is uniformly bounded in~$A$, so that $p$ is uniformly bounded
below, by~\Ref{pmin}.

To choose the set~$C_\KKK := \{1,2,\ldots,a_\KKK\}$, note that, from \Ref{uki}--\Ref{Tsi}
and~\Ref{calcs},
\[
   T_{si} \Le \frac{\alp_i}{\alp_s(\bir_s+\dea_s)} \Eq O\Bl
      i^{-1} e^{-\a (i-s)^2/2A} \Br
\]
for all~$i$, and thus
$$
    T_s \Le \frac{1}{\alp_s(\bir_s+\dea_s)}\sum_{i\ge1}\alp_i \Eq O(A^{-1/2}).
$$
Hence, from~\Ref{CK-def}, the choice $\lceil a_\KKK = 2(s+1)\rceil$ corresponds to a value of
$\KKK \le 1$.  For the corresponding value of~$T_\KKK^+$, it is necessary to bound the
expressions for~$T_k$, which is in detail tedious; however, it is not difficult
to deduce that $T_k = O(\log A \vee \log k)$, so that $T_\KKK^+ = O(\log A)$.
From Theorem~\ref{quasi-longtime}, it now follows that $\law_s(X(t))$ is close
to~$\pds$ for all times~$t$ such that
\[
     A^2 (\log A)^3 \ \ll\ t \ \ll\ A^{-1/2}e^{A/2\a}.
\]
Furthermore, $\m(T) = O(\log A)$ for all return distributions concentrated on sets
of the form $\{1,2,\ldots,A^m\}$ for any fixed exponent~$m$, and Theorem~\ref{main-TV}
thus shows that the corresponding equilibrium distributions~$\pi^\m$ are all exponentially
close to~$\pds$ as $A\to\infty$ --- indeed, $\m$ would have to have extraordinarily
long tails for anything else to be the case.  Hence the fact that this process has
infinitely many \qsd s should not be interpreted as showing any kind of practical
instability, at least for large~$A$: there is an effective long time stable distribution,
and it is extremely close to~$\pds$.

Rather similar analyses could be undertaken for a variety of other well-known
models.  An analogue of the Beverton \& Holt~(1957) model would have $\bbb(x)
=  b/(1 + x/m)$ for $b>d$ and $m>0$, that of Hassell~(1975) would have
$\bbb(x) = b/(1 + x/m)^c$, and that of Maynard--Smith \& Slatkin~(1973) would
have $\bbb(x) = b/(1 + (x/m)^c)$.  The qualitative conclusions would be
entirely similar.

\section{Markov population processes}\label{MPP}
\setcounter{equation}{0}
In this section, we consider Markov population
processes $X_N := (X_N(t),\,t\ge0)$, $N\ge1$, taking values in~$\integ_+^d$,
for some $d\ge1$. In many applications,
the components represent the numbers of individuals of a particular type or species,
with a total of~$d$ types possible.  The process evolves
as a Markov process with state-dependent transitions
\eq\label{1.0}
    X \ \to\ X + J \quad\mbox{at rate}\quad N\a_J(N^{-1}X),\qquad X \in \integ_+^d,\ J \in \jj,
\en
where $\jj \subset \integ^d$ is a fixed finite set, and we define $J_* := \max_{j\in\jj}|J|$.
Density dependence is reflected in the fact that the arguments of the functions~$\a_J$
are counts normalised by the `typical size'~$N$.
The functions~$\a_J\colon\, \integ_+  \to \re_+$ are assumed to be twice continuously
differentiable on~$\re_+^d$, and to be such as to ensure that~$X_N$ is {\it locally irreducible\/};
that is, the number of steps required to get from any state~$X\neq0$ to any
of its lattice neighbours $X+e\uj$, $1\le j\le d$, with positive probability,
is uniformly bounded.

Such processes satisfy a law of large numbers (Kurtz, 1970), expressed in terms of the
system of {\it deterministic equations\/}
\eq\label{determ}
   \frac{d\x}{dt} \Eq \sjj J\a_J(\x) \ =:\ F(\x),\qquad \x\in\re^d;
\en
here, $\x(t)$ approximates $x_N(t) := N^{-1}X_N(t)$, and the quantity~$F$ represents
the infinitesimal average drift of the components of the random process.
We now suppose that $F(c) = 0$ for some $c\in\re^d$ with $c_j > 0$, $1\le j\le d$,
and that all the eigenvalues of the matrix of derivatives~$DF(c)=:A$ have negative
real parts.  In
this case, $c$ is a locally stable equilibrium of the deterministic system~\Ref{determ},
and, if~$X_N$ is started with $N^{-1}X_N(0)$ close to~$c$, the law of large numbers
implies that $x_N(t)$ remains close to~$c$, in the sense that
\eq\label{LLN-equilib}
  \sup_{0\le t\le T}|x_N(t) - c| \ \to_d\ 0,
\en
for any finite $T > 0$.  The central limit theorem in Kurtz (1971) also shows
that
\eq\label{CLT-equilib}
  N^{1/2}(x_N(\cdot) - c) \ \Rightarrow\ x \qquad \mbox{in}\quad D[0,T],
\en for any $T>0$, where~$x$ is a Gaussian process whose stationary
distribution has zero mean and covariance matrix~$\S$ satisfying
\eq\label{Sigma-eqn}
     A\S + \S A^T + \s^2(c) \Eq 0,
\en
where $\s^2(x) := \sjj JJ^T \a_J(x)$.  Here, we complement this approximation,
by using Theorem~\ref{quasi-longtime-2} to show that the distribution of~$X_N(t)$ is
close in total variation, for time periods that become extremely long as~$N$ increases,
to the equilibrium distribution~$\tpds_{N}$ of a truncated
process~$\tX_N$, which is returned to a specified state~$s := s_N$ near~$Nc$ whenever it
leaves a neighbourhood~$C'(N)$ of~$Nc$.  Of course, this distribution, appropriately
centred and normalized, converges to $MVN_d(0,\S)$ as $N\to\infty$.

In order to prove such a result, we need to define the neighbourhood~$C'(N)$,
and to show that the quantities $(1-\tr_{C'(N)})$, $1/\tT_{C'(N)}$ and $\tT_{C'(N)}^+$
appearing in Theorem~\ref{quasi-longtime-2} can be suitably bounded. The inequality
\eq\label{first-bnd}
    1/\tT_{C'(N)} \Le q_{s_N} \Eq N\sjj \a_J(N^{-1}s_N)
\en
is immediate.  For the remaining bounds, we use
Lyapounov--Foster--Tweedie methods (Meyn \& Tweedie, 1993).
We write $y := N^{1/2}(x-c)$ and $r_y^2 := y^T Vy$, where the positive
definite symmetric matrix~$V$ is to be chosen later, and we first consider the process
$y_N(\cdot) := N^{1/2}(x_N(\cdot)-c)$
stopped when $|r_{y_N(t)}| \ge c_0 N^{1/2}$, for~$c_0$ also to be chosen later.
Then, for~$y$ such that $|r_y| < c_0 N^{1/2}$, the generator~$\aaa$ of the
Markov process acting on a real function~$g(y)$ takes the form
\eq\label{aa-def}
  (\aaa g)(y) \Eq \sjj N\a_J(c+N^{-1/2}y)\{g(y + N^{-1/2}J) - g(y)\}.
\en
Using Taylor's expansion on~$g$, for $|r| < c_0 N^{1/2}$, we have
\eqa
   &&\left|(\aaa g)(y)
      - \sjj N\a_J(c+N^{-1/2}y)\{N^{-1/2}J^TDg(y) + \half N^{-1}D^2g(y)[J\ut]\}\right|\non\\
   &&\quad\Le N^{-1/2}\h_3(y;g),                \label{1st-bnd}
\ena
where
\[
    \h_3(y;g) \Def \sjj \a_J^{*0} J_*^3  \,\sup_{|u|\le N^{-1/2}J_*}\|D^3g(y+u)\|,
\]
and $\a_J^{*0} := \sup_{x\colon (x-c)^TV(x-c) \le c_0^2} \a_J(x)$.
Similarly, expanding~$\a_J$, we obtain
\eqa
   \phantom{H}
    \left|N^{1/2}\sjj \a_J(c+N^{-1/2}y) J^TDg(y) - y^TA^TDg(y)\right| &\le& N^{-1/2}\h_1(y;g),                \label{2nd-bnd}
\ena
where we have used the facts that $\sjj J\a_J(c) = F(c) = 0$ and that
$\sjj D\a_J(c)J^T = A^T$, and where
\[
    \h_1(y;g) \Def \|Dg(y)\|J_* \sjj \a_J^{*2} |y|^2,
\]
and $\a_J^{*2} := \sup_{x\colon (x-c)^TV(x-c) \le c_0^2} \|D^2\a_J(x)\|$; and then
\eq
   \left|\sjj\{\a_J(c+N^{-1/2}y) - \a_J(c)\} D^2g(y)[J,J]\right| \Le
          N^{-1/2}\h_2(y;g),                \label{3rd-bnd}
\en
with
\[
    \h_2(y;g) \Def \|D^2g(y)\|J_*^2 \sjj \a_J^{*1} |y|,
\]
and $\a_J^{*1} := \sup_{x\colon (x-c)^TV(x-c) \le c_0^2} \|D\a_J(x)\|$.
Thus, if one ignores the error terms, the generator acts on~$g$ as that of
a multivariate Ornstein--Uhlenbeck process,
\eq\label{aa-approx}
   (\aaa g)(y) \ \approx\ y^TA^TDg(y) + \half\, \trace\{\s(c)D^2g(y)\s(c)\},
\en
with drift matrix~$A$ and infinitesimal covariance matrix~$\s^2(c)$.

   We now consider the generator acting on functions~$g$ of the form
$g(y) = F_\e(r)$, where $F_\e(r) := \int_\e^r f(t)\,dt$, and the function~$f$
is non-negative.  This gives
\[
     Dg(y) \Eq f(r)r_y^{-1}Vy;\quad D^2g(y) \Eq r_y^{-1}f(r)V +
       \{f'(r_y) - r_y^{-1}f(r_y)\}r_y^{-2}Vy(Vy)^T.
\]
Thus the first term in the approximation~\Ref{aa-approx} to~$\aaa$ yields
\[
    y^TA^TDg(y) \Eq  f(r_y)r_y^{-1}y^TA^TVy \Eq \half f(r_y)r_y^{-1}y^T\{A^TV + VA\}y.
\]
In order to choose functions~$g$ such that~$g(y_N(t))$ is  a super-martingale,
we would like the right hand side to be negative, which will be the case if~$V$ is chosen
in such a way that the symmetric matrix $(A^TV + VA)$ is negative definite.  One way of
doing so here is to take $V := \S^{-1}$, where~$\S$ is as in~\Ref{Sigma-eqn}, in which
case
\[
   A^TV + VA \Eq \S^{-1}\{\S A^T + A\S\}\S^{-1} \Eq -\S^{-1}\s^2(c)\S^{-1}
\]
is immediately negative definite. The remaining term in~\Ref{aa-approx} then gives
\eqs
  \lefteqn{\half\, \trace\{\s(c)D^2g(y)\s(c)\}} \\
   &&\Eq \half r_y^{-1}f(r_y)\trace\{\s(c)V\s(c)\} \\
     &&\mbox{}\qquad   +  \half\{f'(r_y) - r_y^{-1}f(r_y)\}r_y^{-2}\trace\{\s(c)Vy(Vy)^T\s(c)\} \\
  &&\Eq \half r_y^{-1}f(r_y)\trace\{\ts^2\} + \half\{f'(r_y) - r_y^{-1}f(r_y)\}R(y)\,,
\ens
where $\ts^2 := \S^{-1/2}\s^2(c)\S^{-1/2}$ is positive definite, and
$$
   R(y) \Def \frac{(\S^{-1/2}y)^T \ts^2 \S^{-1/2}y}{y^T\S^{-1}y}
$$
is bounded between its smallest and largest eigenvalues $\g$ and~$\G$.

We begin by taking $f(r) := r^{-m}e^{\b r^2}$, for~$m$ and~$\b$ to be chosen suitably.
Then~\Ref{aa-approx} gives the main part of~$(\aaa g)(y)$ as
\eqs
    \lefteqn{\half e^{\b r_y^2}r_y^{-m+1}\Bigl\{ -r_y^{-2}y^T \S^{-1/2}\ts^2\S^{-1/2} y
      + r_y^{-2}\trace\{\ts^2\}  + \{2\b - r_y^{-2}(m+1)\} R(y) \Bigr\} }\\
    &&\Le -\half e^{\b r_y^2}r_y^{-m+1}\bigl\{ -(\g - 2\b) - r_y^{-2}(m+1-\trace\{\ts^2\})\bigr\}
            \phantom{XXXXXXXX}\\
    &&=:\ -\BBB(r_y),
\ens
say.  We now choose $\b$ and~$m$ in such a way that $2\b < \g$ and $m+1 > \trace\{\ts^2\}$.

For the remainders, we note first, for~$\h_3$, that there exist constants $c_-,c_+$ and~$K$
such that
$$
   \sup_{|u|\le N^{-1/2}J_*}\|D^3g(y+u)\| \Le K \|D^3g(y)\|\quad
        \mbox{for}\ c_-N^{-1/2} \le \|y\| \le c_+ N^{1/2},
$$
and that $\|D^3g(y)\| \le C_3 r_y^{-m-2}\exp(\b r_y^2)\{1 + r_y^4\}$ for some~$C_3$.
Thus, for all~$y$ such that $N^{-1/2} c_3 \le r_y \le c'_3 N^{1/2}$, for
suitable~$c_3,c'_3$, where we also choose $c'_3 \le c_0$,
it follows that $\h_3(y;g) \le \BBB(r_y)/6$.  Similar considerations
for $\h_1(y;g)$ and $\h_2(y;g)$ show that, possibly increasing $c_3$ and decreasing~$c'_3$,
the inequality
\[
    \h_1(y;g) + \h_2(y;g) + \h_3(y;g) \Le \half\BBB(r_y)
\]
holds for all~$y\in B(N^{-1/2}c_3,N^{1/2}c'_3)$, where
\[
   B(\r,R) \Def \{y\colon\, \r \le r_y \le R\}.
\]
Hence, for such~$y$, we always have
\[
    (\aaa g)(y) \Le -\half\BBB(r_y) < 0.
\]
Thus the quantity $F_\e(y_N(t\wedge\htt_{\r,R}))$
is a non-negative super-martingale, for any $N^{-1/2}c_3 \le \r < R \le N^{1/2}c'_3$
and any~$0 < \e \le \r$, where
$$
    \htt_{\r,R} \Def \inf_{t\ge0} \{y_N(t) \notin B(\r,R)\}.
$$

Defining $p(\r,R;r) := \pr[y_N(\htt_{\r,R}) \in B(0,\r)\giv r_{y_N(0)} = r]$, it thus
follows easily from the optional stopping theorem that
\eq\label{p-bnd}
   1 - p(\r,R;r) \Le F_\r(r)/F_\r(R) \Le \frac{4\b R^{m+1}}{(m-1)\r^{m-1}}\,e^{\b(r^2-R^2)},
\en
for $\r,R$ such that $N^{-1/2}c_3 \le \r < R \le N^{1/2}c'_3$ and $2\b R(R-\r) \ge 1$,
with the last condition ensuring that a simple lower bound for $F_\r(R)$ is valid.
So take
\eq\label{C(N)-def}
   C'(N) \Def \{X\colon N^{-1/2}(X-Nc) \in B(0,N^{1/2}c'_3)\},
\en
and let $s_N\in\integ_+^d$ be the closest lattice point to~$Nc$.  Any path of~$y_N$ starting
in $B(0, N^{-1/2}c_3)$
has positive probability of hitting~$N^{-1/2}(s_N-Nc)$ by taking the most direct path
from $y_N(0)$ to~$N^{-1/2}(s_N-Nc)$, and this probability is uniformly bounded away
from~$0$, by the local irreducibility assumption
on~$X_N$, and because the number of possible values of~$y_N(\cdot)$ in $B(0, N^{-1/2}c_3)$
is uniformly bounded as
$N$ varies.  Hence
\eq\label{t-hat-bnd}
    P[\t_{\{s_N\}} < \htt_{0,2N^{-1/2}c_3} \giv y_N(0) \in B(0, N^{-1/2}c_3)] \ >\ \d
\en
for some $\d > 0$.  If the complementary event occurs, then $B(0, N^{-1/2}c_3)$
is hit again by~$y_N$ before it leaves $B(0,N^{1/2}c'_3)$ with probability at least
\[
    1 - {K N^m}\,e^{-(c'_3)^2 \b N },
\]
for some~$K$, in view of~\Ref{p-bnd}.  It thus follows that
\eq\label{r_N-bnd}
   1-\tr_{C'(N)} \Le \d^{-1}{K N^m}\,e^{-(c'_3)^2 \b N } \Eq O(e^{-\b' N}),
\en
for any  $0 < \b' < (c'_3)^2\b$.

In order to control the mean time to hitting~$s$ for the process~$\tX_N$,
we take $f(r) := r^{-m} + \th r$, for $m$ large enough and~$\th$ small enough positive.
Then $(\aaa g)(y)$ once again has two principal negative contributions, the first,
bounded above by $-\half\g\th r_y^2$, coming from the drift term, and the second,
bounded above by $-\half mr_y^{-(m+1)}$, from the variance term.  The former dominates
all positive terms for $r_y \ge r_0$, for some fixed~$r_0$,
and the second then dominates for the smaller
values of~$r_y$, if~$m$ is chosen large enough; the quantities $\h_l(y;g)$, $1\le l\le 3$,
are treated as before, and the upper and lower bounds for~$r_y$ can be left unchanged.
Hence, in the same range of~$y$, we always have
\[
    (\aaa g)(y) \Le -\d' < 0,
\]
for some $\d' > 0$.  Applying the optional stopping theorem then yields
\[
    \d' \ex\{\htt_{\r,R} \giv r_{y_N(0)} = r\} \Le F_\r(r)
      \Le \frac{\r^{-m+1}}{m-1} + \half\th r^2 ,
\]
if $m > 1$, uniformly in $\r,R$ such that $N^{-1/2}c_3 \le \r < R \le N^{1/2}c'_3$.
Take the extreme values for $\r$ and~$R$. Then since, for this~$R$,
the process~$\tX_N$ is returned directly to~$s_N$ if
$y_N(\htt_{\r,R}) \notin B(0,R)$, and since the mean time to either hitting
$N^{-1/2}(s_N - Nc)$ or leaving $B(0,2N^{-1/2}c_3)$, starting within $B(0,N^{-1/2}c_3)$,
is uniformly bounded by some $c_1 < \infty$, a regenerative argument much as above
shows that
\[
   \tT_{C'(N)}^+ \Le c_0\{N^{(m-1)/2} + N\} + c_1 + (1-\d)\tT_{C'(N)}^+,
\]
for~$\d$ as in~\Ref{t-hat-bnd}, and hence that, uniformly in~$N$,
\eq\label{T-tilde-bnd}
    \tT_{C'(N)}^+ \Le C\{N^{(m-1)/2} + N\},
\en
for some $C < \infty$.

Collecting the above bounds, we have enough to prove the following theorem.

\begin{theorem}\label{MPP-thm}
   Suppose that~$X$ is a Markov population process with transition rates~$N\a_J$ as given
in~\Ref{1.0}, and that the~$\a_J$ are such as to ensure that~$X_N$ is locally irreducible.
Suppose also that $F(c) = 0$ for some $c\in\re^d$ with $c_j > 0$, $1\le j\le d$,
and that all the eigenvalues of the matrix of derivatives~$DF(c)$ have negative
real parts.  Then there exist $\a$, $\b_1$, $\b_2$ and $c'_3>0$, and $C_1$, $C_2$
and~$C_3 < \infty$, depending only on the parameters of the process and not on~$N$,
such that, for all~$t$,
\[
   \dtv(\tpds_{C'(N)},\law_s(X(t))) \Le C_1 t Ne^{-\b_1 N}
       + C_2 t^{-1/2}N^{1+3\a/2} + C_3 e^{-\b_2 t N^{-\a}},
\]
where $s=s_N$ is the nearest lattice point to~$Nc$, and $C'(N)$ is as defined
in~\Ref{C(N)-def}.
\end{theorem}

\proof
All that is needed is to apply the estimate given in Theorem~\ref{quasi-longtime-2}.
An upper bound on
$(1 - \tr_{C'(N)})$ is given in~\Ref{r_N-bnd}; a bound on~$\tT_{C'(N)}^+$ is given
in~\Ref{T-tilde-bnd}; and $1/\tT_{C'(N)}$ is bounded in~\Ref{first-bnd}.  For the
exponent~$\b_1$, any~$\b'$ as for~\Ref{r_N-bnd} can be taken; $\a = \max\{(m-1)/2,1\}$
as in~\Ref{T-tilde-bnd}; and~$\b_2$ can be taken to be $(1 - \log 2)/\{32C\}$, for~$C$
as in~\Ref{T-tilde-bnd}.
\ep

In view of Theorem~\ref{MPP-thm}, the equilibrium
distribution $\tpds_{N}$ is a very good approximation in total
variation to $\law_s(X_N(t))$, provided that~$t$ is bounded below by a
suitable power of~$N$ and above by a quantity growing exponentially with~$N$.

The lower bound given here for the time at which the quasi-equilibrium approximation
becomes accurate is very pessimistic.  The main reason is that the general
coupling strategy used to prove Theorems \ref{quasi-longtime} and~\ref{quasi-longtime-2}
can be very inefficient in specific instances, and is so here.  Better results
could be expected by using the methods to be found in Roberts \& Rosenthal (1996).

\section*{Acknowledgement}
ADB wishes to thank the School of Mathematical Sciences at Monash University,
the Australian Research Council Centre of Excellence
for Mathematics and Statistics of Complex Systems, and the Institute for Mathematical
Sciences at the National University of Singapore, for their hospitality and financial
support while part of this work was accomplished.


\begin{thebibliography}{3}
\ignore{
\bibitem{And91}
{\sc W.\ J.\ Anderson} (1991),
{\em Continuous--time {M}arkov chains: an appli\-cations\---oriented approach\/}.
Springer-Verlag, New York.
}

\bibitem{BP10}
{\sc A.\ D.\ Barbour \& P.\ K.\ Pollett} (2010)
Total variation approximation for quasi-stationary distributions.
{\em J.\ Appl.\ Probab.\/} {\bf 47}, 934--946.

\ignore{
\bibitem{a03}
{\sc F.\ Arrigoni} (2003), Deterministic approximation of a stochastic
metapopulation model. {\em Adv.\ Appl.\ Prob.\/} {\bf 35},  691--720.
%
\bibitem{b76}
{\sc A.\ D.\ Barbour} (1976), Quasi-equilibrium distributions in Markov population processes.
{\it Adv. Appl. Prob.\/} {\bf 8},  296--314.
%
\bibitem{b80}
{\sc A.\ D.\ Barbour} (1980), Equilibrium distributions for Markov population processes.
{\it Adv. Appl. Prob.\/} {\bf 12}, 591--614.
}
\ignore{
\bibitem{bs09}
{\sc A.\ D.\ Barbour \& S.\ N.\ Socoll} (2009),
Local limit approximations for Markov Population Processes.
{\it J.\ Appl.\ Probab.\/} {\bf  46}, 690--708.
%
\bibitem{bs10}
{\sc A.\ D.\ Barbour \& S.\ N.\ Socoll} (2010),
Translated Poisson approximation to equilibrium distributions of Markov
population processes.
{\it Meth.\ Comp.\ Appl.\ Probab.\/} (to appear).
%
\bibitem{Bar60}
{\sc M.\ S.\ Bartlett} (1960),
{\em Stochastic population models in ecology and epidemiology\/},
Methuen, London.
}

\bibitem{BH57}
{\sc R.\ J.\ H.\ Beverton \& S.\ J.\  Holt} (1957),
On the Dynamics of Exploited Fish Populations.
{\em Fishery Investigations\/} Series II Volume XIX, Ministry of Agriculture,
Fisheries and Food.

\bibitem{DS65}
{\sc J.\ N.\ Darroch \& E.\ Seneta} (1965),
On quasi-stationary distributions in absorbing discrete-time Markov chains.
{\em J.\ Appl.\ Probab.\/} {\bf 2}, 88--100.

\ignore{
\bibitem{Ewe63}
{\sc W.\ J.\ Ewens} (1963),
The diffusion equation and a pseudo-distribution in genetics.
{\em J.\ Roy.\ Statist.\ Soc., Ser~B\/} {\bf 25}, 405--412.
%
\bibitem{Ewe64}
{\sc W.\ J.\ Ewens} (1964),
The pseudo-transient distribution and its uses in genetics.
{\em J.\ Appl.\ Probab.\/} {\bf 1}, 141--156.
%
\bibitem{FKMP95}
{\sc P.\ Ferrari, H.\ Kesten, S.\ Mart\'{i}nez \& P.\ Picco} (1995),
Existence of quasi-stationary distributions. A renewal dynamic approach.
{\em Ann.\ Probab.\/} {\bf 23}, 501--521.
}

\ignore{
\bibitem{HK95}
{\sc K.\ Hamza and F.\ C.\ Klebaner} (1995), Conditions for integrability of Markov chains.
{\it J. Appl. Prob.\/} {\bf 32}, 541--547.
%
\bibitem{TK81}
{\sc T.\ G.\ Kurtz} (1981), {\it Approximation of population processes.\/}
CBMS-NSF Regional Conf.\ Series in Appl.\ Math.\ {\bf 36}, SIAM, Philadelphia.
}

\bibitem{Hass75}
{\sc M.\ P.\ Hassell} (1975)
Density--dependence in single--species populations.
{\em J.\ Anim.\ Ecol.\/} {\bf 45}, 283--296.

\bibitem{Kur70}
{\sc T.\ G.\ Kurtz} (1970)
Solutions of ordinary differential equations as limits of pure jump Markov processes.
{\em J.\ Appl.\ Probab.\/} {\bf 7}, 49--58.

\bibitem{Kur71}
{\sc T.\ G.\ Kurtz} (1971)
Limit theorems for sequences of jump Markov processes approximating ordinary
differential processes.
{\em J.\ Appl.\ Probab.\/} {\bf 8}, 344--356.

\ignore{
\bibitem{l02}
{\sc T.\ Lindvall} (2002), {\em Lectures on the coupling method}, 2nd Edn. Dover.
}

\bibitem{MSS73}
{\sc J.\ Maynard-Smith \& M.\ Slatkin} (1973)
The stability of predator--prey systems.
{\em Ecology\/} {\bf 54}, 384--391.

\bibitem{MT93}
{\sc S.\ P.\ Meyn \& R.\ L.\ Tweedie} (1993)
Stability of Markovian processes III: Foster--Lyapunov criteria for
continuous time processes.
{\em Adv.\ Appl.\ Probab.\/} {\bf 25}, 518--548.

\ignore{
\bibitem{Pi74}
{\sc J.\ W.\ Pitman} (1974),
Uniform rates of convergence for Markov chain transition probabilities.
{\em Z.\ Wahr\-schein\-lich\-keits\-theorie verw.\ Geb.\/} {\bf 29}, 193--227.
}

\bibitem{RR96}
{\sc G.\ O.\ Roberts \& J.\ S.\ Rosenthal} (1996)
Quantitative bounds for convergence rates of continuous time Markov processes.
{\em Electr.\ J.\ Probab.\/} {\bf 1}, Paper no.\ 9.

\bibitem{SVJ66}
{\sc E.\ Seneta \& D.\ Vere--Jones} (1966)
On quasi-stationary distributions in discrete--time Markov chains with a denumerable
infinity of states.
{\em J.\ Appl.\ Probab.\/} {\bf 3},  403--434.

\bibitem{vD91}
{\sc E.\ A.\ Van Doorn} (1991),
Quasi-stationary distributions and convergence to
quasi-stationarity of birth-death processes.
{\em Adv.\ Appl.\ Probab.\/} {\bf 23}, 683--700.

\bibitem{Ve38}
{\sc P.-F.\ Verhulst} (1838),
Notice sur la loi que la population poursuit dans son accroissement.
{\it Correspondance math\'ematique et physique\/} {\bf 10}, 113--121.

\bibitem{Ya47}
{\sc A.\ M.\ Yaglom} (1947),
Certain limit theorems of the theory of branching processes.
{\it Doklady Akad.\ Nauk SSSR (N.S.)\/} {\bf 56}, 795--798.

\end{thebibliography}
\enddocument
\end